# CASTELNUOVO-MUMFORD REGULARITY: EXAMPLES OF CURVES AND SURFACES

MARC CHARDIN AND CLARE D'CRUZ*

## INTRODUCTION

The behaviour of Castelnuovo-Mumford regularity under "geometric" transformations is not well understood. In this paper we are concerned with examples which will shed some light on certain questions concerning this behaviour. One simple question which was open (see e.g. [R]) is: May the regularity increase if we pass to the radical or remove embedded primes? By examples, we show that this happens. As a by-product we are also able to answer some related questions. In particular, we provide examples of licci ideals related to monomial curves in $\mathbf{P}^3$ (resp. in $\mathbf{P}^4$) such that the regularity of their radical is essentially the square (resp. the cube) of that of the ideal.

It is well known that the regularity cannot increase when points (embedded or not) are removed. Hence, in order to construct examples where on removing an embedded component the regularity increases, we have to consider surfaces. More surprisingly, we find an irreducible surface such that, after embedding a line into it, the depth of the coordinate ring increases!

Another important concept to understand is the limit of validity of Kodaira type vanishing theorems. The Castelnuovo-Mumford regularity of the canonical module of a reduced curve is 2. An analogous result holds true for higher dimensional varieties with isolated singularities (in characteristic zero), thanks to Kodaira vanishing. As a consequence one can give bounds for Castelnuovo-Mumford regularity (see [CU]). In [Mum], Mumford proves that for an ample line bundle $\mathcal{L}$ on a normal surface $\mathcal{S}$, $H^1(\mathcal{S}, \mathcal{O}_\mathcal{S} \otimes \mathcal{L}^{-1}) = 0$. He remarks that this is false if $\mathcal{S}$ doesn't satisfy $S_2$ (*i.e.* $\mathcal{S}$ is not Cohen-Macaulay) and asks if the $S_2$ condition is sufficient. The first counter-example was given in [AJ]. In this article, we show that counter-examples satisfing $S_1$ give rise to counter-examples satisfying $S_2$. We then provide monomial surfaces whose canonical module has large Castelnuovo-Mumford regularity, so that this vanishing fails. We give a simple proof to show that if $\mathcal{S}$ satisfies $R_1$, then $H^1(\mathcal{S}, \omega_\mathcal{S} \otimes \mathcal{L}) = 0$

*Date*: October 17, 2002.

2000 *Mathematics Subject Classification*. 13D45, 13C40, 13D02, 13D40.

* Supported by Ministère de la recherche.

and that $H^1(\mathcal{S}, \mathcal{O}_\mathcal{S} \otimes \mathcal{L}^i)$ has constant dimension for $i < 0$ if in addition $\mathcal{L}$ is very ample; this dimension measures the defect of Macaulayness of $\mathcal{S}$.

These monomial surfaces also provide examples of complete intersection surfaces with two reduced irreducible components such that one of them has bigger regularity than that of the complete intersection. This doesn't happen in dimension one.

## 1. Motivations and general setting

One of our motivations was from the following result that the first author learned from a joint reading of [L] with D. Cutkosky.

**Theorem 1.1.** [L, 6.5]. *Let $I$ be an homogeneous ideal of codimension $r$ in $R := k[X_0, \ldots, X_n]$, where $k$ is a field of characteristic 0. We denote by $I^{top}$ the intersection of primary components of $I$ of minimal codimension. Assume that $I$ is generated in degree at most $d$. Then there exists an ideal $J$ of codimension $r$, such that*

(1) $J \subseteq I^{top}$,
(2) $I^r \subseteq J$,
(3) *for any prime $\mathfrak{p}$ containing $I^{top}$ such that $(R/I)_\mathfrak{p}$ is regular, $I_\mathfrak{p} = J_\mathfrak{p}$,*
(4) $\operatorname{reg}(J) \leq r(d-1) + 1$.

*Proof.* Following [L], set $\mathfrak{a} := \tilde{I}$ and $J := \oplus_\mu H^0(\mathbf{P}_k^n, \mathcal{J}(\mathfrak{a}^r)(\mu))$. Then (1) follows from [L, 10.23], (2) from [L, 10.16], (3) from [L, 6.8] and (4) from [L, 4.10 (ii)].

As Lazarsfeld remarks, this result implies that if $X := \operatorname{Proj}(R/I)$ coincides with a smooth equidimensional scheme $Y$ outside finitely many points, then $\operatorname{reg}(I_Y) \leq r(d-1)+1$ (this result generalizes [BEL] also follows from [CU], which also treats certain cases of one dimensional singularities). This remark is an immediate consequence of the following well-known lemma applied to $R/J$,

**Lemma 1.2.** *Let $M$ be a finitely generated graded $R$-module and $\mathfrak{p}$ an homogeneous prime of $R$ such that $\dim R/\mathfrak{p} \leq 1$. Then,*

$$\operatorname{reg}(M/H_\mathfrak{p}^0(M)) \leq \operatorname{reg} M.$$

*Proof.* Let $N := M/H_\mathfrak{p}^0(M)$. Note first that $H_\mathfrak{p}^0(N) = 0$ because $\mathfrak{p} \subset \mathfrak{m}$. As $\dim H_\mathfrak{p}^0(M) \leq 1$, we get isomorphisms $H_\mathfrak{m}^i(M) \simeq H_\mathfrak{m}^i(N)$ for $i \geq 2$, and a surjective map $H_\mathfrak{m}^1(M) \to H_\mathfrak{m}^1(N)$. $\qquad\square$

In remark [L, 6.9] the following question arises: Is $\operatorname{reg}(I^{top}) \leq \operatorname{reg}(J)$ if $I$ defines an irreducible scheme ? With some more generality, one may ask the following questions:

Q1: Does the Lemma extend to primes of positive dimension ?

Q2: Is it possible that the regularity increases if one removes embedded components of positive dimension ?

Q3: Is it possible that the regularity increases if one replaces the ideal by its radical ?

Q4: Is it possible that the depth decreases if one replaces the ideal by its radical ?

Q5: Is it possible that the depth decreases if we one remove embedded components of positive dimension ?

The answer to some of these questions are well-known. It is easy to see that answers to Q1 and Q4 are "no" and "yes", respectively, and it is not difficult to cook-up counter-examples.

The examples in section 2 and section 3 show that the answer to Q2, Q3 and Q5 is "yes", even if we impose the condition that the corresponding scheme is irreducible. Our examples will be mainly curves and surfaces. In particular we get the following result,

**Proposition 1.3.** *There exists an irreducible and reduced surface* $\mathcal{S} \subset \mathbf{P}^5$ *and another scheme* $\mathcal{S}'$ *supported on* $\mathcal{S}$ *which coincides with* $\mathcal{S}$ *outside one line and such that:*
*1)* $\operatorname{reg}(\mathcal{S}') < \operatorname{reg}(\mathcal{S})$,
*2)* $\operatorname{depth}(R/I_{\mathcal{S}'}) > \operatorname{depth}(R/I_{\mathcal{S}})$.

## 2. Counter-examples built from curves.

*In this section $k$ will be a field of any characteristic.*

The monomial curve parametrized by $x_0 = 1, x_1 = w^{a_1}, \ldots, x_n = w^{a_n}$ on the affine chart $X_0 = 1$ of $\mathbf{P}^n_k$ will be called "the monomial curve of degrees $(a_1, \ldots, a_n)$".

For question Q3, there are simple counter-examples for curves.

**Example 2.1.** Let $I := (y^2 z - x^2 t, z^4 - x t^3) \subset k[x, y, z, t]$. Then $I$ defines a complete intersection curve of $\mathbf{P}^3_k$ which has three components: the monomial curve $(1, 6, 8)$, the reduced line $x = z = 0$, and a triple line supported on $z = t = 0$. This can be checked since the ideal of the residuals of each of the two lines is the ideal of the $2 \times 2$ minors of the matrices

$$\begin{pmatrix} z^3 & t^3 \\ y^2 & xt \\ x & z \end{pmatrix} \qquad \eta := \begin{pmatrix} xt^2 & z^3 \\ x^2 & y^2 \\ z & t \end{pmatrix}.$$

In particular the ideal $J := I : (z, t) = I + (xy^2t^2 - x^2z^3)$ and has the following resolution:

$$0 \longrightarrow R[-6]^2 \xrightarrow{\eta} R[-5] \oplus R[-4] \oplus R[-3] \longrightarrow J$$

so that the regularity of $J$ is 5 (compared to 6 for $I$). Also $\sqrt{J} = \sqrt{I}$.

The radical of $J$ is $K := J + (x^4z^2 - xy^4t)$, which has regularity 6 and the following minimal free resolution:

$$0 \longrightarrow R[-8] \xrightarrow{\phi} R[-7]^2 \oplus R[-6]^2 \xrightarrow{\psi} R[-6] \oplus R[-5] \oplus R[-4] \oplus R[-3] \longrightarrow K$$

where $\phi$ and $\psi$ are given by the matrices

$$\phi := \begin{pmatrix} -y^2 \\ -x^2 \\ t \\ z \end{pmatrix} \qquad \psi := \begin{pmatrix} xt^2 & -z^3 & xy^2t & -x^2z^2 \\ -x^2 & y^2 & 0 & 0 \\ z & -t & -x^2 & y^2 \\ 0 & 0 & z & -t \end{pmatrix}.$$

This can be checked by using the Buchsbaum-Eisenbud criterion: after verifying that $\psi \circ \phi = 0$, it suffices to find two $3 \times 3$ minors of $\psi$ without common factor to know that this sequence resolves $K$. Now $K$ is saturated (the resolution has length 3) and is Cohen-Macaulay on the punctured spectrum since $I_1(\phi)$ is $(x, y, z, t)$-primary. Also $(x^4z^2 - xy^4t)^2 \in J$ so that $\sqrt{K} = \sqrt{J}$ and $K$ is unmixed. It remains to compute the degree of $K$. The resolution of $K$ gives the Hilbert-Poincaré series of $R/K$

$$\frac{1 - t^3 - t^4 - t^5 + t^6 + 2t^7 - t^8}{(1-t)^4} = \frac{10}{(1-t)^2} - \frac{19}{(1-t)} + Q(t)$$

which completes the argument.                                    $\square$

This example also provides an answer to case Q2 as follows,

**Example 2.2.** Let $\mathfrak{b}$ be the ideal of the monomial curve $(1, 6, 8)$. It is easy to check along the same lines as in the above example that

$$0 \longrightarrow R[-8] \oplus R[-7] \xrightarrow{\phi} R[-7]^2 \oplus R[-6]^3 \oplus R[-5]$$

$$\xrightarrow{\psi} R[-6] \oplus R[-5] \oplus R[-4]^2 \oplus R[-3] \longrightarrow \mathfrak{b}$$

where $\phi$ and $\psi$ are given by the matrices

$$\phi := \begin{pmatrix} y^2 & 0 \\ t & y^2 \\ z & x^2 \\ -x & 0 \\ 0 & t \\ 0 & -z \end{pmatrix} \qquad \psi := \begin{pmatrix} t^2 & -y^2t & xz^2 & z^3 & y^4 & x^3z \\ z & x^2 & -y^2 & xt & 0 & 0 \\ x & 0 & 0 & y^2 & 0 & 0 \\ 0 & z & -t & 0 & x^2 & y^2 \\ 0 & 0 & 0 & 0 & z & t \end{pmatrix}.$$

is the minimal free resolution of $\mathfrak{b}$ (see [Br] or [BH] for the general case of monomial curves in $\mathbf{P}^3$). In particular, the regularity of this curve is 6.

The ideal $zJ$ has regularity $\operatorname{reg}(J) + 1 = 6$ and its radical is $(z) \cap \mathfrak{b}$. Now $z$ doesn't divide zero in $R/\mathfrak{b}$ so that $(z) \cap \mathfrak{b} = (z)\mathfrak{b}$ and therefore $\operatorname{reg}((z) \cap \mathfrak{b}) = \operatorname{reg}(\mathfrak{b}) + 1 = 7 > \operatorname{reg}(zJ)$. Note that $\sqrt{zJ}$ differs from $zJ$ only at the embedded primes. $\qquad\square$

**Example 2.3.** Consider the family of ideals $I_{m,n} := (x^mt - y^mz, z^{n+2} - xt^{n+1}) \subset k[x,y,z,t]$. Then,

**Lemma 2.4.** *For $m, n \geq 1$, $\operatorname{reg}(I_{m,n}) = m + n + 2$ and $\operatorname{reg}(\sqrt{I_{m,n}}) = mn + 2$.*

*Proof.* The first equality is clear since $I_{m,n}$ is a complete intersection. For the second, we claim that the complex below is a graded free resolution of $\sqrt{I_{m,n}}$,

$$0 \longrightarrow \bigoplus_{i=2}^{n} R[-mi - n + i - 4] \xrightarrow{\phi_{mn}} \bigoplus_{i=1}^{n} R[-mi - n + i - 3]^2$$

$$\xrightarrow{\psi_{mn}} R[-m-1]\bigoplus_{i=0}^{n} R[-mi - n + i - 2] \xrightarrow{\gamma_{mn}} \sqrt{I_{m,n}} \ ,$$

where $\gamma_{mn} := (y^mz - x^mt \ \ f_0 \ \cdots \ f_n)$ and $\phi_{mn}, \ \psi_{mn}$ are given by the matrices

$$\phi_{mn} := \begin{pmatrix} C_m & 0 & \cdots & 0 \\ C & C_m & \ddots & \vdots \\ 0 & C & \ddots & 0 \\ \vdots & \ddots & \ddots & C_m \\ 0 & \cdots & 0 & C \end{pmatrix} \quad \psi_{mn} := \begin{pmatrix} B_{1mn} & B_{2mn} & \cdots & B_{nmn} \\ L_m & 0 & \cdots & 0 \\ L & L_m & \ddots & \vdots \\ 0 & L & \ddots & 0 \\ \vdots & \ddots & \ddots & L_m \\ 0 & \cdots & 0 & L \end{pmatrix}$$

and

$$C := \begin{pmatrix} t \\ z \end{pmatrix}, \quad C_m := \begin{pmatrix} -y^m \\ -x^m \end{pmatrix}, \quad L := \begin{pmatrix} z & -t \end{pmatrix}, \quad L_m := \begin{pmatrix} -x^m & y^m \end{pmatrix},$$

$$B_{imn} := \begin{pmatrix} xy^{(i-1)m}t^{n-i+1} & -x^{(i-1)m}z^{n+2-i} \end{pmatrix}$$

$$f_i := x^{im}z^{n+2-i} - xy^{im}t^{n-i+1}$$

Now $LC = L_m C_m = L C_m + L_m C = B_{imn} C_m + B_{(i+1)mn} C = 0$ so that $\phi_{mn} \circ \psi_{mn} = 0$. Also

$$(y^m z - x^m t)B_{imn} + f_{i-1}L_m + f_i L = (0 \ 0)$$

so that it is a complex. The Buchsbaum-Eisenbud criterion shows the exactness of this complex, the saturation and the unmixedness of the image of $\gamma_{mn}$, noticing that the ideal of maximal minors of the first map is $(x, y, z, t)$-primary and the $(n + 1) \times (n + 1)$ of the second map contains $z^n(z^{n+2} - xt^{n+1})$ and $t^{n-1}(x^m t - zy^m)$. Clearly $x^m t - zy^m$ and the $f_i$'s belong to the ideal of the monomial curve $\mathcal{C}$ of degrees $(1, m(n + 1), m(n + 2))$ and to the ideals $(x, z)$ and $(z, t)$. From the resolution it is easy to see that the degree of the image of $\gamma_{mn}$ is $m(n + 2) + 2$. This implies that $\mathrm{Im}(\gamma_{mn}) = I_{\mathcal{C}} \cap (x, z) \cap (z, t)$.    □

**Example 2.5.** Let $I_{\mathcal{C}}$ be the ideal of the monomial curve of degree $(1, mn^2, mn(n+1), m(n + 1)^2)$ in $\mathbf{P}^4 = \mathrm{Proj}(k[x, y, z, u, v])$ with $m \geq 1$ and $n \geq 3$. Then $I_{\mathcal{C}}$ contains the complete intersection $\mathfrak{I} := (y^m u^2 - x^m zv, z^{n+1} - xu^n, u^{n+1} - xv^n)$. Note that $\mathfrak{I} \subset (x, u^2, z^{n+1})$ and that $J_{m,n} := \mathfrak{I} : (x, z^{n+1}, u^2) = \mathfrak{I} + (y^m v^n - x^{m-1}zu^{n-1}v)$ is not contained in $(x, z, u)$.

Now $\mathrm{reg}(J_{m,n}) = \mathrm{reg}(\mathfrak{I}) - 1 = m + 2n + 1$ and

(1) $J_{m,n} = I_{\mathcal{C}} \cap K$, where $K$ is the $(z, u, v)$-primary component of $\mathfrak{I}$,
(2) $M_1 := y^{mn} - x^{mn^2-1}z$ and $M_2 := y^{m(n^2-2n-1)}v - x^{m(n^2-2n-1)-1}z^2$ are minimal generators of $I_{\mathcal{C}}$,
(3) $\mathrm{reg}(I_{\mathcal{C}}) = mn^2$ and $\mathrm{reg}(\sqrt{J_{m,n}}) = m(n^2 - 2n - 1) + 1$.

*Proof.*    (1) We first compute the multiplicity of $K$. On the affine chart $x = 1$, after inverting $y$ we have $\tilde{\mathfrak{I}} = (au^2 - zv, z^{n+1} - u^n, u^{n+1} - v^n)$. If $z$, $u$ or $v$ is invertible and $a \neq 0$, then the multiplicity of this ideal is one. We need to compute the multiplicity $\mu$ of this ideal at the origin over $k(a)$. Let us consider the family of complete intersection ideals $I_{a,b} := (au^2 - bzv, z^{n+1} - wu^n, u^{n+1} - wv^n)$ in $k[z, u, v, w]$ over $\mathbf{P}^1_k = \mathrm{Proj}(k[a, b])$. It can be easily verified that if $(a, b) \neq (1 : 0)$ or $(0 : 1)$, then $I_{a,b}$ consists of a $(z, u, v)$-primary component, a component in $w = 0$ and a point $p_{a,b} = (a^{-n^2}b^{n+n^2} : b^n : a^n : a^{2n+1}b^{-n-1})$. Hence, the multiplicity at the origin is constant and in fact is $\mu$ for $ab \neq 0$, is $\mu$ for $(a,b) = (1:0)$ and is $\mu - 1$ when $(a,b) = (0:1)$. Now $I_{1,0} := (u^2, z^{n+1} - wu^n, u^{n+1} - wv^n) = (u^2, z^{n+1}, wv^n)$ which has multiplicity $2n(n+1)$ at the origin. Therefore $\mu = 2n(n+1)$.

The degree of $J_{m,n}$ is $(m+2)(n+1)^2 - 2(n+1) = m(n+1)^2 + \mu$ which is the degree of $I_{\mathcal{C}} \cap K \supset J_{m,n}$. This implies the asserted equality.

(2) Recall that ideals of monomial curves are minimally generated by binomials. It is easy to see that the only binomial of the form $y^\beta - z^\gamma u^\delta v^\epsilon$ in the ideal $\tilde{I}_{\mathcal{C}}$ of the affine monomial curve of degree at most $mn^2$ is $y^{mn^2} - z$. Therefore $M_1$ is not in the ideal $((I_{\mathcal{C}})_{<mn^2})$ which is the homogenization of $((\tilde{I}_{\mathcal{C}})_{<mn^2})$. Assume that $M_2$ is not a minimal generator. Then we can write:

$$M_2 = y^{m(n^2-2n-1)}v - x^{m(n^2-2n-1)-1}z^2 = \sum_i \pm Q_i(N_i - P_i)$$

where $N_i - P_i$ is a minimal binomial generators of $I_{\mathcal{C}}$ of degree at most $m(n^2 - 2n - 1)$ and $Q_i$ is a monomial of complementary degree. Either $N_i$ or $P_i$ has to be of the form $x^\alpha z$ or $x^\alpha z^2$. We claim that the only binomials of the form $z - y^\beta u^\delta v^\epsilon$ or $z^2 - y^\beta u^\delta v^\epsilon$ in the ideal $\tilde{I}_{\mathcal{C}}$ of degree at most $mn^2$ are $\tilde{M}_1 = z - y^{mn^2}$, $\tilde{M}_2 = z^2 - y^{m(n^2-2n-1)}v$ and $\tilde{M}_3 = z^2 - y^{m(n^2-n)}u$. This can be seen from the equations $mn^2 = \beta + mn^2\delta + m(n+1)^2\epsilon$ or $2mn^2 = \beta + mn^2\delta + m(n+1)^2\epsilon$, respectively. The assertion follows from the fact that the degrees of $\tilde{M}_1$ and $\tilde{M}_3$ are bigger than the one of $M_2$.

(3) It follows from (2) that $I_{\mathcal{C}} + (z, u, v) = (y^{mn^2}, z, u, v)$, so that setting $R := k[x, y, z, u, v]$ and $L := (z, u, v)$ we have an exact sequence,

$$0 \to R/I_{\mathcal{C}} \cap L \to R/I_{\mathcal{C}} \oplus R/L \to k[x,y]/(y^{mn^2}) \to 0.$$

Therefore, the Hilbert series of $R/I_{\mathcal{C}}$ and $R/I_{\mathcal{C}} \cap L$ may be computed one from another. Namely,

$$H_{R/I_{\mathcal{C}} \cap L}(t) = H_{R/I_{\mathcal{C}}}(t) + \frac{1}{(1-t)^2} - \frac{1 - t^{mn^2}}{(1-t)^2}.$$

Now the natural epimorphisms $R/\mathfrak{J} \to R/I_{\mathcal{C}}$ and $R/\mathfrak{J} \to R/I_{\mathcal{C}} \cap L$ show that the second local cohomology modules of $R/I_{\mathcal{C}}$ and $R/I_{\mathcal{C}} \cap L$ vanishes in degrees $\geq a(R/\mathfrak{J}) = m + 2n - 1$. As (2) implies that $\mathrm{reg}(I_{\mathcal{C}}) \geq mn^2$ and $\mathrm{reg}(\sqrt{J_{m,n}}) \geq m(n^2 - 2n - 1) + 1$, the regularity of these ideals may be computed from the Hilbert function of $R/I_{\mathcal{C}}$.

In the appendix, we compute this function in degrees $\geq mn + m + 2n$, from which (3) follows. □

**Remark 2.6.** The ideals $I_{m,n}$ of Example 2.3 or $J_{m,n}$ of Example 2.5 are such that $\mathrm{reg}(zI_{m,n}) = m + n + 3$ and $\mathrm{reg}(zJ_{m,n}) = m + 2n + 2$. They differ from $(z) \cap I_{\mathcal{C}} = (z)I_{\mathcal{C}}$ –for $\mathcal{C}$ equal respectively to the monomial curves $(1, m(n+1), m(n+2))$ and $(1, mn^2, mn(n+1), m(n+1)^2)$– only along lines embedded in $z = 0$. [Alternatively, the ideals $z(I_{m,n} : (xt, z))$ or $z(J_{m,n} : (u, v, z))$ can be used to decrease the regularity by one and at the same time keep the same property.]

Now $\mathrm{reg}((z) \cap I_{\mathcal{C}}) = \mathrm{reg}(I_{\mathcal{C}}) + 1 = mn + 3$ (resp. $\mathrm{reg}((z) \cap I_{\mathcal{C}}) = mn^2 + 1$), which follows from the estimates of the regularity of the corresponding monomial curves.

## 3. Canonical module of surfaces

Let $X \subset \mathbf{P}_k^n$ be a projective scheme of codimension $r$ over a field $k$. We set $\omega_X := \mathcal{E}xt^r_{\mathcal{O}_{\mathbf{P}^n}}(\mathcal{O}_X, \omega_{\mathbf{P}^n})$.

**Proposition 3.1.** *Let $\mathcal{S}$ be a projective unmixed scheme over a field. Let $\pi : \mathcal{S}' = \mathrm{Spec}\,\mathcal{E}nd(\omega_{\mathcal{S}}) \to \mathcal{S}$ be its $S_2$-ification and $\mathcal{L}$ be an ample line bundle on $\mathcal{S}$. Then $\mathcal{S}'$ is projective, $\pi^*\mathcal{L}$ is ample and for any $i$,*

$$H^i(\mathcal{S}, \pi^*\mathcal{L} \otimes \omega_{\mathcal{S}'}) = H^i(\mathcal{S}, \mathcal{L} \otimes \omega_{\mathcal{S}}).$$

*Proof.* As $\pi$ is finite, $\mathcal{S}'$ is projective and $\pi^*\mathcal{L}$ is ample. Notice that $\pi_*\omega_{\mathcal{S}'} = \omega_{\mathcal{S}}$ so that by the finiteness of $\pi$ and the projection formula we have,

$$H^i(\mathcal{S}, \pi^*\mathcal{L} \otimes \omega_{\mathcal{S}'}) = H^i(\mathcal{S}, \pi_*(\pi^*\mathcal{L} \otimes \omega_{\mathcal{S}'})) = H^i(\mathcal{S}, \mathcal{L} \otimes \pi_*\omega_{\mathcal{S}'}) = H^i(\mathcal{S}, \mathcal{L} \otimes \omega_{\mathcal{S}}).$$

$\square$

**Corollary 3.2.** *Let $\mathcal{S}$ be a surface defined over a field of characteristic zero. If $\mathcal{S}$ is regular in codimension one, then for any ample line bundle $\mathcal{L}$, $H^1(\mathcal{S}, \mathcal{L} \otimes \omega_{\mathcal{S}}) = 0$.*

*Proof.* As $\mathcal{S}'$ is the normalization of $\mathcal{S}$ in this case, this follows from the Proposition 3.1 and [M]. $\square$

**Remark 3.3.** *Let $\mathcal{S}$ be an unmixed scheme of dimension $d$ over a field of characteristic zero. Note that more generally one has $H^{d-j}(\mathcal{S}, \mathcal{L} \otimes \omega_{\mathcal{S}}) = 0$ for $j \le \min\{i, d-1\}$ if $\mathcal{S}$ satisfies $R_i$ and $\mathcal{L}$ is ample (see e.g. [CU, 1.3]). On the other hand, if $\mathcal{S}$ satisfies $S_{i+1}$, then $H^{d-j}(\mathcal{S}, \mathcal{L} \otimes \omega_{\mathcal{S}}) \simeq H^j(\mathcal{S}, \mathcal{L}^{-1})$ for $j < i$ and there is a surjective map $H^{d-i}(\mathcal{S}, \mathcal{L} \otimes \omega_{\mathcal{S}}) \to H^i(\mathcal{S}, \mathcal{L}^{-1})$.*

*Therefore, if $\mathcal{S}$ satisfies $R_i$ and $S_{i+1}$, then $H^j(\mathcal{S}, \mathcal{L}^{-1}) = 0$ for $\mathcal{L}$ ample and $j \leq \min\{i, d-1\}$. This is proved in* [AJ].

*Note also that these vanishings extend to a big and nef line bundle on a scheme having at most rational singularities locally in codimension $i$ for the first and satisfying $S_{i+1}$ in addition for the second.*

In [M], Mumford asks if the vanishing he proves extends to varieties satisfying $S_2$, in particular to irreducible Cohen-Macaulay surfaces. The first counter-example was given in [AJ]. Also Proposition 3.1 shows that if this where true, then it would hold for irreducible surfaces in general. The following examples show that $H^1(\mathcal{S}, \mathcal{L} \otimes \omega_{\mathcal{S}}) \neq 0$ for some monomial surfaces in $\mathbf{P}^5$ with $\mathcal{L} := \mathcal{O}_{\mathbf{P}^5}(1)$.

**Example 3.4.** *Consider the surface parametrized in $\mathbf{A}^5$ by*

$$x_1 = a^5, \ x_2 = b^6, \ x_3 = a^4b, \ x_4 = ab^2, \ x_5 = a^2b^5$$

*with $k = \mathbf{Q}$ or $\mathbf{Z}/101\mathbf{Z}$. Then the homogenisation of the ideal of this irreducible surface gives an homogeneous prime $\mathfrak{p}$ in $R := k[X_0, \ldots, X_5]$ defining its closure $\mathcal{S} \subset \mathbf{P}^5$. Using the "Macaulay 2" software, we get that $H^2_{\mathfrak{m}}(\omega_{R/\mathfrak{p}})$ has four socle elements of respective degrees $-1, -1, 0, 1$. In particular we have $\dim_k H^1(\mathcal{S}, \omega_{\mathcal{S}}(1)) = 1$.*

<div align="right">□</div>

Certain binomial ideals corresponding to monomial surfaces of higher degrees provide more counter-examples.

**Example 3.5.** Consider the surface in $\mathbf{A}^5$ parametrized by

$$x_1 = a^{12}, \ x_2 = b^8, \ x_3 = ab^7, \ x_4 = a^5b, \ x_5 = a^9b^4$$

with $k = \mathbf{Q}$ or $\mathbf{Z}/101\mathbf{Z}$, and let $\mathfrak{q}$ be the ideal of its projective closure in $R$. Then $H^2_{\mathfrak{m}}(\omega_{R/\mathfrak{q}})$ has a socle element of degree 5, and $\mathrm{reg}(\omega_{R/\mathfrak{q}}) = 7$. Also $\mathfrak{q}$ has regularity 32 and $I := (\mathfrak{q}_{\leq 21})$ has regularity 24 and is of the form $I = \mathfrak{q} \cap J$ where $J$ is $(X_1, X_2, X_3, X_5)$-primary and embedded in $I$. More surprising is the fact that $\mathrm{depth}(R/\mathfrak{q}) = 1$ and $\mathrm{depth}(R/I) = 2$.

<div align="right">□</div>

**Remark 3.6.** The characteristics 0 and 101 that we choose for computing in "Macaulay 2" are of course arbitrary. We didn't check that these counter-examples work in any characteristic, which seems to require hard checking for no significant improvement.

**Remark 3.7.** By [CU, 4.2], any link $Y$ of $X := \mathrm{Proj}(R/\mathfrak{q})$ by a complete intersection $Z$ of degrees $d_1, d_2, d_3$ satisfy

$$d_1 + d_2 + d_3 = \mathrm{reg}(Y) > \mathrm{reg}(Z) = d_1 + d_2 + d_3 - 3.$$

This link may be chosen reduced and irreducible if the $d_i$'s are $\geq 21$. In contrast, note that in the case of curves $\mathrm{reg}(Y) < \mathrm{reg}(X \cup Y)$ if $X$ is reduced and $X \cup Y$ is a complete intersection.

**Proposition 3.8.** *Let $\mathcal{S}$ be a surface in $\mathbf{P}^n$ and $B$ be its homogeneous coordinate ring. Let $\omega_B := \mathrm{Ext}_R^{n-2}(B, \omega_R)$, $N := \mathrm{Ext}_R^{n-1}(B, \omega_R)$, $\mathcal{N} := \widetilde{N}$ and $\mathcal{S}' := \mathrm{Proj}(\mathrm{End}(\omega_B))$ be its Macaulayfication. Then for $\mu \in \mathbf{Z}$*

(1) $H_\mathfrak{m}^2(\omega_B) \cong H_\mathfrak{m}^0(N)$;

(2) *There is an exact sequence,*

$$0 \longrightarrow H^0(\mathcal{S}, \mathcal{O}_\mathcal{S}(\mu)) \longrightarrow H^0(\mathcal{S}, \mathcal{O}_{\mathcal{S}'}(\mu)) \longrightarrow H^1(\mathcal{S}, \mathcal{N}(-\mu)) \longrightarrow 0 \;;$$

(3) $H^1(\mathcal{S}, \omega_\mathcal{S}(-\mu)) \simeq H^1(\mathcal{S}, \omega_{\mathcal{S}'}(-\mu))$, $H^0(\mathcal{S}, \mathcal{N}(-\mu)) \simeq k^\ell$, *where $\ell$ is zero if and only if $\mathcal{S}$ is Cohen-Macaulay and equal to $\deg N$ if $\mathcal{S}$ is not Cohen-Macaulay. And there is an exact sequence*

$$0 \longrightarrow H^1(\mathcal{S}, \omega_\mathcal{S}(-\mu)) \longrightarrow H^1(\mathcal{S}, \mathcal{O}_\mathcal{S}(\mu))^\vee \longrightarrow H^0(\mathcal{S}, \mathcal{N}(-\mu))$$
$$\longrightarrow H^1(\mathcal{S}, \mathcal{N}(-\mu)) \longrightarrow 0 \;;$$

(4) *Assume that $\mathcal{S}$ is reduced. Then*
   (a) $H_\mathfrak{m}^2(B)_\mu \simeq k^\ell$, $\forall \mu < 0$ *if and only if $H_\mathfrak{m}^2(\omega_B)_{>0} = 0$;*
   (b) $H_\mathfrak{m}^2(B)_\mu \simeq k^\ell$, $\forall \mu \leq 0$ *if and only if $\mathcal{S}$ is connected in codimension one and $H_\mathfrak{m}^2(\omega_B)_{\geq 0} = 0$.*

*Proof.* Let $L_\bullet$ be a graded free resolution of $\Gamma B := \ker(\mathcal{C}_\mathfrak{m}^1(B) \longrightarrow \mathcal{C}_\mathfrak{m}^2(B)) = \bigoplus_\mu H^0(\mathcal{S}, \mathcal{O}_\mathcal{S}(\mu))$, and set $-^\vee := \mathrm{Hom}(-, \omega_R)$ and $-^* := \mathrm{Homgr}(-, k)$. The double complex $\mathcal{C}_\mathfrak{m}^\bullet(L_\bullet^\vee)$ gives rise to two spectral sequences. They respectively have terms: $'E_1^{pq} = H_\mathfrak{m}^p(L_q^\vee)$ which is zero if $p \neq n + 1$ and else isomorphic to $L_q^*$. Therefore

$$'E_2^{n+1,q} = H^q(H_\mathfrak{m}^{n+1}(L_\bullet^\vee)) \simeq H^q(L_\bullet^*) \simeq H_q(L_\bullet)^*$$

so that $'E_2^{p,q} = 0$ unless $(p, q) = (n + 1, 0)$ and $'E_2^{n+1,0} = (\Gamma B)^*$.

On the other hand, $''E_2^{p,q} = H_\mathfrak{m}^p(\mathrm{Ext}_R^q(\Gamma B, \omega_R))$ which is zero unless $q \geq n-2$ and $p \leq n + 1 - q$. As $\mathrm{Ext}_R^q(\Gamma B, \omega_R) = 0$ for $q \geq n$ and $\mathrm{Ext}_R^q(\Gamma B, \omega_R) = \mathrm{Ext}_R^q(B, \omega_R)$ for $q \leq n - 1$ there are only four non-zero modules $''E_2^{p,q}$ for $(p, q) = (0, n-1), (1, n-1), (2, n-2), (3, n-2)$ and two non trivial maps:

$$''d_2^{2,n-2} : H_\mathfrak{m}^2(\omega_B) \longrightarrow H_\mathfrak{m}^0(N) \qquad ''d_2^{3,n-2} : H_\mathfrak{m}^3(\omega_B) \longrightarrow H_\mathfrak{m}^1(N)$$

where the first one is an isomorphism and the second one has a kernel isomorphic to $(\Gamma B)^*$.

This proves (1) and (2). For (3), $H^1(\mathcal{S}, \omega_{\mathcal{S}}(-\mu)) \simeq H^1(\mathcal{S}, \omega_{\mathcal{S}'}(-\mu))$ by Proposition 3.1 and the second claim follows from the fact that $N$ is of dimension at most one. Also the exact sequence follows from (1), (2), the isomorphism $N \simeq H^2_{\mathfrak{m}}(B)^*$ and the standard exact sequence

$$0 \longrightarrow H^0_{\mathfrak{m}}(N)_{-\mu} \longrightarrow N_{-\mu} \longrightarrow H^0(\mathcal{S}, \mathcal{N}(-\mu)) \longrightarrow H^1_{\mathfrak{m}}(N)_{-\mu} \longrightarrow 0 \ .$$

If $\mathcal{S}$ is reduced, so is $\mathcal{S}'$ and (2) implies that $H^1_{\mathfrak{m}}(N)_\mu = 0$ for $\mu > 0$ and also $H^1_{\mathfrak{m}}(N)_0 = 0$ if and only if $H^0(\mathcal{S}, \mathcal{O}_{\mathcal{S}}) = H^0(\mathcal{S}, \mathcal{O}_{\mathcal{S}'})$ which in turn is equivalent to $\mathcal{S}$ being connected in codimension one. Therefore (4) follows from the exact sequence of (3). $\qquad \square$

**Corollary 3.9.** *Let $\mathcal{S}$ be an embedded projective surface with isolated singularities. Then $H^1(\mathcal{S}, \mathcal{O}_{\mathcal{S}}(\mu))$ has constant dimension for $\mu < 0$.*

And we have seen counter-examples with non isolated singularities.

# 4. Appendix

The Hilbert function of the monomial curve $(1, mn^2, mn(n+1), m(n+1)^2)$ is determined by the number of monomials in the ideal

$$I^\alpha \ = \ (x^{m(n+1)^2}, x^{m(n+1)^2-1}y, x^{m(2n+1)}y^{mn^2}, x^{m(n+1)}y^{mn(n+1)}, y^{m(n+1)^2})^\alpha.$$

We are interested in high degrees. Since it is enough to consider the exponents of $y$, we will work on the corresponding subsets in $\mathbf{N}$. If $E$ and $F$ are two such subsets, $E + F := \{e + f \mid e \in E, \ f \in F\}$, also $kE := E + \cdots + E$ ($k$ times), and $[a, b] := \{n \in \mathbf{N} \mid a \leq n \leq b\}$. We will write:

$$I(\alpha) = \alpha\{0, 1, mn^2, mn(n+1), m(n+1)^2\} = \bigcup_{i=0}^{\alpha} mn^2(\alpha - i) + J(\alpha - i, i)$$

where

$$J(\beta, i) := \beta\{0, mn, m(2n+1)\} + [0, i].$$

In what follows we put $\beta := \alpha - i$ and $I_i(\alpha) := mn^2(\alpha - i) + J(\alpha - i, i)$.

$\bullet$ For $i \geq m(n+1) - 1$, $J(\beta, i) = [0, (2n+1)m\beta + i]$, and the first and last elements of $I_i(\alpha)$ are

$$
\begin{aligned}
F(\beta, i) &:= \beta mn^2 \\
L(\beta, i) &:= \beta m(n+1)^2 + i.
\end{aligned}
$$

Now,

$$L(\beta, i) - F(\beta + 1, i - 1) + 1 = \alpha - (\beta - 1) - m(n^2 - 2\beta n - \beta).$$

Hence $I_i(\alpha) \cup I_{i-1}(\alpha)$ is an interval whenever $\beta \geq \frac{mn^2 - \alpha - 1}{m(2n+1) - 1}$. Therefore, if $\alpha \geq mn^2 - 1$,

$$\bigcup_{i \geq m(n+1)-1} I_i(\alpha) = [0, L(\alpha - m(n+1) + 1, m(n+1) - 1)].$$

If $mn^2 - 1 > \alpha \geq mn + m + n$, then it follows that the set $S_\alpha$ of missing exponents of $y$ smaller than $L(\alpha - m(n+1) + 1, m(n+1) - 1)$ has cardinality

$$(a+1)\mu - \binom{a+1}{2}(m(2n+1) - 1),$$

where $\mu := mn^2 - 1 - \alpha$ and $a = \left\lfloor \frac{\mu}{m(2n+1) - 1} \right\rfloor$.

- Let $i < m(n+1) - 1$. Notice that

$$\beta\{0, n, 2n + 1\} = \bigcup_{s=0}^{2\beta} E_s(\beta, n),$$

with

$$E_s(\beta, n) := \begin{cases} \{0, \ldots, \lfloor \frac{s}{2} \rfloor\} + sn & 0 \leq s \leq \beta \\ \{0, \ldots, \lfloor \frac{2\beta - s}{2} \rfloor\} + sn & \beta < s \leq 2\beta \end{cases}$$

which leads us to put $E_s^i(\beta, m, n) := \{me + j \mid e \in E_s(\beta, n), \ 0 \leq j \leq i\}$.

We remark that $\max E_s^i(\beta, m, n)) \geq \min E_{s+1}^i(\beta, m, n) - 1$ if and only if

$$(*) \qquad 2(n - \left\lfloor \frac{i+1}{m} \right\rfloor) \leq s < \beta - 2(n - \left\lfloor \frac{i+1}{m} \right\rfloor)$$

which gives a partition $[0, 2\beta] = a \coprod b \coprod c$, with $a := [0, 2(n - \left\lfloor \frac{i+1}{m} \right\rfloor) - 1]$, $b$ given by the condition $(*)$ and $c := [\beta - 2(n - \left\lfloor \frac{i+1}{m} \right\rfloor), 2\beta]$, and the corresponding partition $J(i, \beta) := A(\beta, i) \coprod B(\beta, i) \coprod C(\beta, i)$.

- For $m - 1 \leq i < m(n+1) - 1$, $B(\beta, i)$ is an interval. We put $F(\beta, i) := \min\{B(\beta, i) + mn^2\beta\}$ and $L(\beta, i) := \max\{B(\beta, i) + mn^2\beta\}$. Now $F(\beta, i)) \leq L(B(\beta - 1, i + 1)) + 1$ if

$$m(\beta - 2n)(2n + 1) \geq m - i - 1$$

which holds true if $\beta \geq 2n$. This is the case if $\alpha \geq mn + m + 2n$. As $L(\alpha - m(n+1) + 1, m(n+1) - 1) \geq F(\alpha - m(n+1) + 2, m(n+1) - 2))$,

$$[0, L(\alpha - m + 1, m - 1)] \cap I(\alpha) = [0, L(\alpha - m + 1, m - 1)] - S_\alpha.$$

- For $0 \leq i \leq m - 1$, we put

$$G_i(\alpha) := \{\gamma \in I_i(\alpha) \mid \gamma > L(\alpha - m + 1, m - 1)\},$$

and notice the following facts,

– $A(\alpha - i, i) + mn^2(\alpha - i) \cap G_i(\alpha) = \emptyset$,

– every element $\gamma \in G_i(\alpha)$ is of the form $\gamma =: km + j$ with $0 \leq j \leq i$ and, if $\alpha \geq mn + m + 2n$, then

$$\gamma \in G_{i-1}(\alpha) \iff j < i.$$

Hence

$$\bigcup_{i=0}^{m-1} G_i(\alpha) = \coprod_{i=0}^{m-1} (H_i(\alpha) := \{\gamma \in G_i(\alpha) \mid \gamma \equiv i \mod m\}).$$

The elements of $E_s^i(\beta, m, n)$ congruent to $i$ mod $m$ are in one to one correspondance with the elements of $E_s(\beta, n)$. All elements of $C(\alpha - i, i) + mn^2(\alpha - i)$ are in $G_i(\alpha)$ and the ones that are in $H_i(\alpha)$ correspond to $s = 2\beta, \ldots, 2\beta - 2n$ and this gives a contribution of $\sum_{k=0}^{2n}(\lfloor \frac{k}{2} \rfloor + 1) = (n+1)^2$ to $H_i(\alpha)$. To count the elements of $B(\alpha - i, i) + mn^2(\alpha - i)$ in $G_i(\alpha)$ one notes that

$$B(\beta, i) + mn^2\beta = \{\gamma \in [F(\beta, i), L(\beta, i)] \mid \exists\, 0 \leq j \leq i, \gamma \equiv j \mod m\},$$

and therefore

$$\begin{aligned}
|H_i(\alpha)| &= (n+1)^2 + \frac{L(\alpha - i, i) - L(\alpha - m + 1, m - 1) + m - 1 - i}{m} \\
&= (n+1)^2 + (n+1)^2(m - 1 - i)
\end{aligned}$$

so that

$$|\bigcup_{i=0}^{m-1} G_i(\alpha)| = (n+1)^2(m + \binom{m}{2}).$$

Finally, if $\alpha \geq mn + m + 2n$, then

$$\begin{aligned}
|I(\alpha)| &= L(\alpha - m + 1, m - 1) + 1 - |S_\alpha| + (n+1)^2(m + \binom{m}{2}) \\
&= m(n+1)^2\alpha - \left[ \binom{m}{2}(n+1)^2 + m(n^2 - 1) \right] - |S_\alpha|
\end{aligned}$$

where $|S_\alpha| = 0$ if $\alpha \geq mn^2 - 1$ and if $\alpha = mn^2 - 1 - \mu$ (with $\mu > 0$):

$$|S_\alpha| = \dim_k H^1_{\mathfrak{m}}(R/I_{\mathcal{C}})_\alpha = (a+1)\mu - \binom{a+1}{2}(m(2n+1) - 1),$$

where $a = \left\lfloor \frac{\mu}{m(2n+1)-1} \right\rfloor$. This implies that the Hartshorne-Rao module has length $\geq m^2 n^5 / 4$.

Institut de Mathématiques, CNRS & Université Paris 6, 4, place Jussieu, F–75252 Paris cedex 05, France

*E-mail address*: chardin@math.jussieu.fr

Chennai Mathematical Institute, 92-G . N.Chetty Road, T. Nagar, Chennai 600 017, India

*Current address*: Institut de Mathématiques, CNRS & Université Paris 6, 4, place Jussieu, F–75252 Paris cedex 05, France

*E-mail address*: clare@cmi.ac.in